\documentclass [11pt] {article}

\usepackage{amsmath,amssymb}

\newtheorem{theorem}{Theorem}

  \newcommand {\qed} {\null \hfill \rule{2mm}{2mm}}

\def\sep{{\rm sep}}

\def \eqalign#1{\null\,\vcenter{\openup\jot\mathsurround 0pt\ialign{\strut
\hfil$\displaystyle{##}$&$\displaystyle{{}##}$\hfil
&&\quad\strut\hfil$\displaystyle{##}$&$\displaystyle{{}##}$\hfil
\crcr#1\crcr}}\,}

\def \ialign {\everycr {} \tabskip 0pt \halign }

\def\Res{{\rm Res}}

\begin {document}

\title{{\Large{\bf Root separation for reducible integer polynomials}}}

\author{Yann Bugeaud and Andrej Dujella}

\date{}
\maketitle

\footnotetext{
 \vspace{1ex} \\
The authors were supported by the French-Croatian bilateral COGITO
project {\it Polynomial root separation}.}

\begin{abstract}
We construct parametric families of (monic) reducible polynomials
having two roots very close to each other.
\end{abstract}

\section{Introduction}

The (na\"\i ve) height $H(P)$ of  an integer
polynomial $P(x)$ is
the maximum of the absolute values of its
coefficients.
For a separable
integer polynomial $P(x)$ of degree $d \ge 2$
and with distinct roots $\alpha_1, \ldots , \alpha_d$,
we set
$$
\sep(P) = \min_{1 \le i < j \le d} \, |\alpha_i - \alpha_j|
$$
and define the quantity $e(P)$ by
$$
\sep(P) = H(P)^{-e(P)}.
$$
Following the notation introduced in \cite{BM2}, for $d \ge 2$, we set
$$
e(d) := \limsup_{{\rm deg}(P) = d, H(P) \to + \infty} e(P)
$$
and
$$
e_{{\rm irr}} (d) := \limsup_{{\rm deg}(P) = d, H(P) \to + \infty} e(P),
$$
where the latter limsup is taken over the irreducible integer polynomials
$P(x)$ of degree $d$.
We further define
$e^*(d)$ and $e_{{\rm irr}}^*(d)$ by restricting to
monic, respectively, monic irreducible integer polynomials, of degree $d$.
Obviously, we have
$$
e(d) \ge e_{{\rm irr}}(d) \quad \mbox{and} \quad
e^*(d) \ge e_{{\rm irr}}^*(d), \quad (d \ge 2).
$$
A classical result of Mahler \cite{Mah64}
asserts that $e(d) \le d-1$ for every $d \ge 2$
and it is easy to check that $e_{{\rm irr}}(2)=e(2)=1$
and $e^*(2)=e_{{\rm irr}}^*(2) = 0$.

The determination of the exact values of $e(d), e^*(d), e_{{\rm irr}} (d)$
and $e_{{\rm irr}}^*(d)$ has been investigated by several authors over the
past ten years. Evertse \cite{Ev}
and Sch\"onhage  \cite{Sch} proved, independently, that
$e_{{\rm irr}}(3)=e(3)=2$. Rather surprisingly,
no other known value of $e(d)$ is known.
The following inequalities gather the lower bounds obtained
by Beresnevich, Bernik, Bugeaud, Dujella, G\"otze, and Mignotte in
the four papers \cite{BM1,BM2,BBG,BD}:
$$
e_{{\rm irr}}(d) \geq \frac{d}{2} +\frac{d-2}{4(d-1)}, \quad \mbox{for $d\geq 4$},
$$
$$
e (d)  \geq \frac{d+1}{2}, \quad \mbox{for odd $d\geq 5$},
$$
$$
e^*(3)=e_{{\rm irr}}^*(3) \geq \frac{3}{2}, \quad e^* (5) \ge 2,
\quad e_{{\rm irr}}^*(5) \ge \frac{7}{4},
$$
$$
e^*(d)\geq \frac{d}{2}, \quad e_{{\rm irr}}^*(d) \geq \frac{d-1}{2},
\quad \mbox{for even $d\geq 4$},
$$
and
$$
e^*(d)\geq \frac{d-1}{2}, \quad
e_{{\rm irr}}^*(d)\geq \frac{d}{2} +\frac{d-2}{4(d-1)} - 1, \quad \mbox{for odd $d\geq 7$}.
$$
See also \cite{DP} for a constructive proof that $e^*(4)\geq 2$
and \cite{Dub} for the study of analogous quantities defined in terms
of the Remak height (instead of the na\"\i ve height).

The aim of the present paper is to improve all known lower bounds for $e(d)$, $e^*(d)$ and
$e_{{\rm irr}}^* (d)$ for $d$ sufficiently large.

\begin{theorem} \label{tm1}
For any integer $d\geq 4$, we have
$$
e(d)  \geq \frac{2d-1}{3}.
$$
\end{theorem}

We obtain a slightly weaker lower bound when we restrict our
attention to monic polynomials.

\begin{theorem} \label{tm2}
For any even positive integer $d\geq 6$,  we have
$$
e^*(d) \geq \frac{2d-3}{3}.
$$
For any odd positive integer $d\geq 7$, we have
$$
e^*(d) \geq \frac{2d-5}{3}.
$$
\end{theorem}

Roughly speaking, all the previously known
lower bounds were of order $d/2$.
There are many other questions on integer polynomials of degree $d$,
or on algebraic numbers of degree $d$, for which the answer is known to lie
somewhere between $d/2$ and $d$. The most celebrated one is the problem
of Wirsing on the approximation to transcendental real numbers by algebraic
numbers of degree at most $d$; see Chapter 3 of \cite{BuLiv}.
We stress that the lower bounds in Theorems \ref{tm1} and \ref{tm2}
are of the order $2d/3$. As far as we are aware, this is the first time where
an estimate of order $\theta d$ with $\theta > 1/2$ is obtained for such kind of questions.

As pointed out in \cite{Bu1}, irreducible polynomials with close roots are useful
to investigate the difference between Mahler's and Koksma's classifications of real numbers.
However, it does not seem to us that Theorem \ref{tm1} could be applied to this question
to improve Corollary 1 of \cite{BD}.

To prove Theorems 1 and 2 we construct parametric families
of integer polynomials. For an integer $d \ge 4$,
the reducible polynomials arising in the proof
of Theorem 1 are products of a linear polynomial
$L_n (x) = (n^2 + 3n +1)x - (n+2)$ with an irreducible
polynomial $p_{d-1, n} (x)$ of degree $d-1$ and height
of order $n$. We then show that $p_{d-1, n} (x)$ has a root
$y_{d,n}$ very close to the root $x_n = (n+2)/(n^2 + 3n + 1)$
of $L_n(x)$. Say differently, we construct a parametric family
$(y_{d,n})_{n \ge 1}$ of algebraic numbers of degree $d-1$ which are very
well approximated by a rational number with large height (in comparison to the
height of $y_{d,n}$). This means that it is then possible to apply Bombieri's version
of the Thue--Siegel principle \cite[Theorem 4]{Bom}
to the anchor pair $(y_{d,n}, x_n)$ to derive a rather good effective irrationality measure
for $y_{d,n}$ when $n$ is sufficiently large in comparison to $d$.
In his paper, Bombieri used the polynomials $x^d - n x  + 1$,
which have a root very close to the rational number $1/n$.

In our last result, we improve the known lower bound
for $e_{{\rm irr}}^* (d)$ when $d$ is large enough.

\begin{theorem} \label{tm3}
For any positive integer $d\geq 4$,  we have 
$$
e_{{\rm irr}}^* (d) \geq \frac{d}{2}-\frac{1}{4}.
$$
\end{theorem}

Throughout the next sections, the constants implied by the symbols $O$, $\ll$ and
$\gg$ can be explicitly computed, are
independent of the parameter $n$, and depend at most on the degree $d$.

\section{Reducible polynomials: Proof of Theorem \ref{tm1}}

We want to construct a one-parametric sequence of integer
polynomials $p_{d,n}(x)$ of degree $d$ having a root very close to the rational number
$x_n=(n+2)/(n^2+3n+1)$. Then the polynomials
$$
P_{d,n}(x)=((n^2+3n+1)x-(n+2)) p_{d-1,n}(x)
$$
will have two roots very close to each other.
We define the sequence $p_{d,n}(x)$ recursively by
$$
\eqalign{
& p_{0,n} (x) = -1, \quad p_{1,n} (x) = (n+1)x-1,  \cr
& p_{d,n} (x) = (1+x) p_{d-1,n}(x) +x^2 p_{d-2,n} (x). \cr} \eqno (1)
$$

We claim that
$$
p_{d,n}\left(\frac{n+2}{n^2+3n+1}\right) = \frac{(-1)^{d-1}}{(n^2+3n+1)^d}.    \eqno (2)
$$
Indeed, (2) is clearly true for $d=0$ and $d=1$. Assume now that $d \ge 1$
is an integer for which (2)
holds for $p_{d-1,n} (x)$ and $p_{d,n} (x)$.
Then we deduce from the recursion (1) that
$$
\eqalign{
p_{d,n}\left(\frac{n+2}{n^2+3n+1}\right)
& = \frac{(-1)^{d-2}}{(n^2+3n+1)^{d-1}} \cdot \frac{n^2+ 4n+3}{n^2+3n+1}  \cr
&   \, \, \, \, \, \, \, \, +
\frac{(-1)^{d-3}}{(n^2+3n+1)^{d-2}} \cdot \frac{n^2+ 4n+4}{(n^2+3n+1)^2} \cr
& = \frac{(-1)^{d-1}}{(n^2+3n+1)^d}, \cr}
$$
as claimed.

We now show that for sufficiently large $n$ the polynomial $p_{d,n} (x)$ has a root
between $x_n$ and
$$
z_{d,n}=x_n + \frac{(-1)^{d}}{n(n^2+3n+1)^d}.
$$
Observe that
$$
(-1)^{d-1} p_{d,n}(x_n) = \frac{1}{(n^2+3n+1)^d} > 0.
$$
By Rolle's theorem, there exists $z'_{d,n}$ between $x_n$ and $z_{d,n}$
such that
$$
p_{d,n}(z_{d,n})=p_{d,n}(x_n)+(z_{d,n}-x_n)p'_{d,n}(z'_{d,n}).
$$
It follows easily by induction that
$$
p_{d,n}(x)=-1+(n-d+2)x+((d-1)n-(d-1)(d-2)/2)x^2+\cdots
$$
Since $x_n=1/n+O(1/n^2)$, we have $p'_{d,n}(z'_{d,n}) = n + d + O(1/n)$.
Thus, for sufficiently large $n$, we get $p'_{d,n}(z'_{d,n}) > n$.
This implies that
$$
(-1)^{d-1}p_{d,n}(z_{d,n}) = \frac{1}{(n^2+3n+1)^d} - \frac{1}{n(n^2+3n+1)^d}p'_{d,n}(z'_{d,n}) < 0.
$$
Therefore, the polynomial $P_{d,n}(x)=((n^2+3n+1)x-(n+2)) p_{d-1,n}(x)$ has two close roots:
$x_n$ and $y_{d,n}$, which is between $x_n$ and $z_{d-1,n}$.
This yields
$$
{\rm sep}(P_{d,n}) \leq |x_n - y_{d,n} | \leq
\frac{1}{n(n^2+3n+1)^{d-1}} \leq \frac{1}{n^{2d-1}},   \eqno (3)
$$
when $n$ is large enough. Since the height of $P_{d,n} (x)$
is bounded from above by $n^3$ times a number depending only on $d$, this gives
$$
e(d)\geq \frac{2d-1}{3},
$$
by letting $n$ tend to infinity.
The proof of Theorem \ref{tm1} is complete.
\qed

\medskip

By Gelfond's inequality (see e.g.  \cite[Lemma A.3]{BuLiv}),
the height of $y_{d,n}$ is $\ll  n$.
Liouville's inequality (see e.g.  \cite[Theorem A.1]{BuLiv}) then implies that
$$
|x_n - y_{d,n}| \gg  n^{-2 \delta-1},
$$
where $\delta$ denotes the degree of $y_{d,n}$.
Combined with (3), this gives $\delta = d-1$ and establishes that
the polynomial $p_{d-1,n} (x)$ must be irreducible.
Also, (3) shows that Liouville's inequality
$$
|x_n - y_{d,n}| \gg  n^{-2d+1}
$$
is sharp in terms of $n$.


\section{Reducible monic polynomials: Proof of Theorem \ref{tm2}}

In order to get a family of monic polynomials with similar
separation properties as the family $P_{d,n}(x)$, we replace the
linear non-monic polynomial $L_n(x)=(n^2+3n+1)x-(n+2)$ by
the monic quadratic polynomial
$$
K_n(x)=x^2 - (n^2+3n+1)x+(n+2).
$$
Thus, we want to construct a one-parametric sequence of integer
polynomials $q_{d,n}(x)$ of degree $d$ having a root very close
to the root $y_n=1/n+O(1/n^2)$ of $K_n(x)$.
Then the polynomials
$$
Q_{d,n}(x)=(x^2 - (n^2+3n+1)x+(n+2)) q_{d-2,n}(x)
$$
will have two roots very close to each other.

For $d \ge 0$ even, we define the sequence $q_{d,n}(x)$ recursively by
$$
\eqalign{
& q_{0,n} (x) = 1, \quad q_{2,n} (x) = x^2-(n+1)x+1, \cr
& q_{d,n} (x) = (2x^2 + x + 1) q_{d-2,n} (x) -  x^4 q_{d-4,n} (x). \cr} \eqno (4)
$$
We claim that $q_{d,n} (x)- q_{d-2,n} (x) q_{2,n} (x)$ is divisible by $K_n (x)$.
This is easy to check for $d=2$ and $d=4$,
and then the claim follows by induction using the recursion (4).
This yields that
$$
q_{d,n}(y_n) = q_{d-2,n}(y_n)q_{2,n}(y_n) = (q_{2,n}(y_n))^{d/2},
$$
for $d \ge 2$ even. From
$$
y_n=1/n-1/n^2+2/n^3-4/n^4+8/n^5+ O(1/n^6),
$$
we get $q_{2,n}(y_n)=1/n^4 + O(1/n^5)$ and hence
$$
q_{d,n}(y_n) = 1/n^{2d} + O(1/n^{2d+1}).
$$

We now show that for sufficiently large $n$ the polynomial $q_{d,n} (x)$ has a root
between $y_n$ and $w_{d,n}=y_n + \frac{2}{n^{2d+1}}$.
We have
$$
q_{d,n}(y_n) = 1/n^{2d} + O(1/n^{2d+1}) > 0.
$$
By Rolle's theorem, there exists $w'_{d,n}$ between $y_n$ and $w_{d,n}$ such that
$$
q_{d,n}(w_{d,n})=q_{d,n}(y_n)+(w_{d,n}-y_n)q'_{d,n}(w'_{d,n}).
$$
It follows easily by induction that
$$
q_{d,n}(x)=1+(-n+d/2-2)x+((-d/2+1)n+(d^2-2d+8)/8)x^2+\cdots
$$
Since $y_n=1/n+O(1/n^2)$, we have $q'_{d,n}(w'_{d,n}) = -n - d/2 + O(1/n)$.
Thus, for sufficiently large $n$, we get $q'_{d,n}(w'_{d,n}) < -n$.
This implies
$$
q_{d,n}(w_{d,n}) = 1/n^{2d} + O(1/n^{2d+1}) + \frac{2}{n^{2d+1}}q'_{d,n}(w'_{d,n}) < 0.
$$

Thus, the polynomial $Q_{d,n}(x)=(x^2 - (n^2+3n+1)x+(n+2)) q_{d-2,n}(x)$
has two close roots:
$y_n$ and $v_{d,n}$, which is between $y_n$ and $w_{d-2,n}$.
This yields
$$
{\rm sep}(Q_{d,n}) \leq \frac{2}{n^{2d-3}},
$$
when $n$ is large enough.
Since $H(Q_{d,n})=O(n^3)$, this gives
$$
e^*(d)\geq \frac{2d-3}{3},
$$
by letting $n$ tend to infinity.

Let now $d$ be odd. Then we define
$$ Q_{d,n}(x)=x(x^2 - (n^2+3n+1)x+(n+2)) q_{d-3,n}(x). $$
This polynomial has two close roots:
$y_n$ and a root lying between $y_n$ and $w_{d-3,n}$.
Thus we get
$$
{\rm sep}(Q_{d,n}) \leq \frac{2}{n^{2d-5}},
$$
for $n$ large enough, and
$$
e^*(d)\geq \frac{2d-5}{3}.
$$
The proof of Theorem \ref{tm2} is complete.
\qed

\section{Irreducible monic polynomials: Proof of Theorem \ref{tm3}}

In this section, we use the polynomials $p_{d,n}(x)$
to construct irreducible monic polynomials
having two very close roots.
Let $F_k$ denote the $k$th Fibonacci number defined by
the recursion $F_0=0$, $F_1=1$ and $F_k=F_{k-1}+F_{k-2}$ for $k \ge 2$.
Note that Fibonacci numbers appear in the asymptotic expansion of
$x_n = (n+2) / (n^2 + 3n + 1)$, namely
$$
x_n = 1/n-1/n^2+2/n^3-5/n^4+13/n^5-34/n^6+ \cdots - (-1)^k F_{2k-3}/n^{k} + \cdots
$$
For $d\geq 0$, we first define
monic polynomials $s_{d,n}(x)$ with a root close to $x_n $ by
$$
s_{d,n} (x)=(-1)^{d-1}(F_{d-1}p_{d,n}(x) - F_d x p_{d-1,n} (x)),
$$
and then monic polynomials with two close roots by
$$
\eqalign{
r_{2d+1,n} (x) & = xs_{d,n}^2(x) + F_d^2 p_{d,n}^2 (x),  \cr
r_{2d,n}(x) & = s_{d,n}^2 (x) + F_{d-1}^2 x p_{d-1,n}^2 (x). \cr} \eqno (5)
$$
We claim that these polynomials are monic.
It suffices to show that this is true for $s_{d,n}(x)$.
Since the leading coefficient of $p_{d,n} (x)$ is $F_d n + F_{d-2}$,
we deduce that the leading coefficient of $s_{d,n}(x)$ is equal to
\begin{eqnarray*}
& (-1)^{d-1} (F_{d-1}(F_d n + F_{d-2})-F_d(F_{d-1}n + F_{d-3})) \\ 
& =
(-1)^{d-1}(F_{d-1}F_{d-2}-F_dF_{d-3}) = 1.
\end{eqnarray*}

From (2) we get
$$
r_{2d+1,n}(x_n )=F_d^2/n^{4d-1} + O(1/n^{4d})
$$
and
$$
r_{2d,n}(x_n )=F_{d-1}^2/n^{4d-3} + O(1/n^{4d-2}),
$$
that is,
$$
r_{d,n}(x_n )=F_{\lfloor (d-1)/2 \rfloor}^2/n^{2d-3} + O(1/n^{2d-2}).
$$

Observe that the degree of the polynomial $r_{d,n} (x)$ is $d$. 
We claim that
$r_{d,n} (x)$ has two complex conjugate roots $v_{d,n}$ and
$\overline{v_{d,n}}$ close to $x_n $,
more precisely they are equal to
$$
\eqalign{
1/n-1/n^2+2/n^3- & 5/n^4+13/n^5- \ldots + \cr
& + (-1)^d F_{2d-5}/n^{d-1} \pm i/n^{(2d-1)/2} + O(1/n^{d}). \cr} \eqno (6)
$$
Indeed, the polynomials $F_d^2 p_{d,n}^2(x)$ and $F_{d-1}^2 x p_{d-1,n}^2(x)$
have double 
roots $y_{d+1,n}$ and $y_{d,n}$, resp., in the disc 
$|x-x_n|\leq 1/n^{2d+1}$,
resp. $|x-x_n|\leq 1/n^{2d-1}$.
Moreover, we have $F_d^2 p_{d,n}^2(x) < xs_{d,n}^2(x)$ when 
$|x-x_n| = 1/n^{2d+1}$
and $F_{d-1}^2 x p_{d-1,n}^2(x) < s_{d,n}^2(x)$ when 
$|x-x_n| = 1/n^{2d-1}$.
Hence, by Rouch\'e's theorem, the polynomial
$r_{d,n} (x)$ has two roots $v$ satisfying $|v-x_n |<1/n^{d-1}$.
Using a complex version of Taylor's theorem,
by writing the roots $v$ of $r_{d,n} (x)$ close to $x_n$
as $v =x_n +\beta+\gamma i$ with $\gamma > 0$, we get 
$$
\eqalign{
0 & =r_{d,n}(x_n +\beta+\gamma i) \cr
& = r_{d,n}(x_n ) + (\beta+\gamma i)r'_{d,n}(x_n )+
(\beta+\gamma i)^2r''_{d,n}(x_n )/2+ O(1/n^{3d-5}). \cr}
$$
Note that we have
$$
r'_{2d+1,n}(x_n ) =\frac{2F_dF_{d+1}}{n^{2d-1}} + O(1/n^{2d}),
$$
$$
r'_{2d,n}(x_n ) =\frac{2F_{d-1}F_{d-2}}{n^{2d-2}} + O(1/n^{2d-1}),
$$
and
$$
r''_{d,n}(x_n ) = 2F_{\lfloor (d-1)/2 \rfloor}^2 n^2 + O(n).
$$
By considering imaginary parts in the above Taylor's formula we get
\begin{eqnarray*}
\beta &=& - \frac{r'_{d,n}(x_n)}{r''_{d,n}(x_n)} + O(\gamma^{-1}r''_{d,n}(x_n)^{-1}n^{-3d+5}) \\
 &=& - \frac{F_{\lfloor (d-1)/2 \rfloor + (-1)^d}}{n^d} + O\left(\frac{1}{n^{d+1}}\right) +
O\left(\frac{1}{\gamma n^{3d-3}}\right).
\end{eqnarray*}
Since the distance between the roots $x_n +\beta+\gamma i$ and  
$x_n +\beta-\gamma i$ of $r_{d,n} (x)$ is equal to $2 \gamma$  
and $H(r_{d,n})=O(n^2)$, 
it follows from Mahler's theorem \cite{Mah64} quoted in the introduction that 
$\gamma \gg 1/n^{2d-2}$. Thus, $\beta = O(1/n^{d-1})$.
Moreover, $\beta = O(1/n^d)$, unless $\gamma \ll 1/n^{2d-3}$.
Looking at real parts leads to
$$
\eqalign{
\frac{F_{\lfloor (d-1)/2 \rfloor}^2}{n^{2d-3}} +
\frac{2\beta F_{\lfloor (d-1)/2 \rfloor}F_{\lfloor (d-1)/2 \rfloor
+ (-1)^d}}{n^{d-2}} + (\beta^2 - \gamma^2) &
F_{\lfloor (d-1)/2 \rfloor}^2 n^2 \cr
& = O\left(\frac{1}{n^{2d-2}}\right).}  \eqno{(7)}
$$
The assumption that $\gamma \ll 1/n^{2d-3}$ leads to a contradiction, 
by considering (7) as a quadratic equation
in $\beta$. Hence, we have $\beta = O(1/n^{d})$. Now (7) gives
$\gamma^2 = 1/n^{2d-1} + O(1/n^{2d})$, i.e.,
$\gamma = 1/n^{(2d-1)/2} + O(1/n^d)$, as claimed in (6). 

Let $R_{d,n}(x)$ be the irreducible factor of $r_{d,n}(x)$ having roots
$v_{d,n}$ and $\overline{v_{d,n}}$. Denote by $\delta$ its degree.
Note that, since $H(r_{d,n})=O(n^2)$, Gelfond's inequality implies that
$$
H(R_{d,n})=O(n^2).    \eqno (8)
$$

Denote by $\Res_{d,n}$ the resultant of
the polynomials $R_{d,n} (x)$ and $L_n(x)$.
Since $\Res_{d,n}$ is a rational integer and $x_n$ is not a root of $R_{d,n} (x)$, we have
$$
|\Res_{d,n}| \ge 1.   \eqno (9)
$$
Furthermore, the definition of the resultant of two polynomials (see e.g. \cite[p. 223]{BuLiv})
implies that (recall that $R_{d,n}(x)$ is monic)
$$
|\Res_{d,n}| \le  (n^2 + 3n+1)^{\delta} \prod_{t : R_{d,n} (t) = 0} \, |t - x_n|.
$$
Since the product of the absolute values of all the roots of $R_{d,n} (x)$
different from $v_{d,n}$ and $\overline{v_{d,n}}$ is bounded
from above by $\sqrt{d+1}$ times
the height of $R_{d,n} (x)$ (see \cite[Exercise A.1]{BuLiv}), we deduce from (6) and (8) that
$$
|\Res_{d,n}| \ll n^{2 \delta} n^{-2d+1} n^2.
$$
Combined with (9), this gives $2 \delta \ge 2d - 3$, thus
$$
\delta \in  \{d, d-1\}.
$$
This implies that either $r_{d,n}(x)$ is irreducible, or it has an integer
root (recall that any rational root of a monic polynomial must be a rational
integer).

Assume that $r_{d,n}(x)$ is reducible, i.e., that it has an integer root.
This integer must divide $r_{d,n}(0)=F_{\lfloor (d-1)/2 \rfloor}^2$,
and, by (5),
it is of the form $-\alpha^2$ for an integer $\alpha$ dividing $F_{\lfloor (d-1)/2 \rfloor}$.
So, for any given degree $d$, there are only finitely many possibilities for the integer root
of $r_{d,n}(x)$.

Write $r_{d,n}(x)$ as
$$
r_{d,n}(x)=A_d(x)n^2+B_d(x)n+C_d(x).
$$
By the definition of $r_{d,n}(x)$ and (1), it follows easily by induction
that the discriminant $B_d(x)^2 - 4A_d(x)C_d(x)$ of $r_{d,n}(x)$ with respect to $n$
is equal to $-4 F_{\lfloor (d-1)/2 \rfloor}^4 x^{2d-1}$.
In particular, this discriminant is nonzero at $x=-\alpha^2$,
and thus at least one of the numbers $A_d(-\alpha^2)$, $B_d(-\alpha^2)$, $C_d(-\alpha^2)$
is nonzero.
This shows that $r_{d,n}(x)$ is irreducible over $\mathbb{Z}[n,x]$.

Furthermore, for given $\alpha$ there are exactly two (rational) values of $n$
such that $-\alpha^2$ is a root of $r_{d,n}(x)$. Hence, for sufficiently large
positive integer $n$, the polynomial $r_{d,n}(x)$ has no integer roots,
and therefore it is irreducible over $\mathbb{Z}[x]$.

Since
$$
{\rm sep}(r_{d,n}) = O(n^{-(d-1/2)}),
$$
we obtain
$$
e_{{\rm irr}}^* (d) \geq \frac{2d-1}{4},
$$
which proves Theorem \ref{tm3}.
\qed

\medskip

\section{Clusters of roots}

Let $d \ge 2$ be an integer.
Mahler's upper bound $e(d) \le d-1$ quoted in the Introduction
is a particular case of the lower bound
$$
\prod_{1\le i<j\le k} \vert \alpha_i -\alpha_j \vert  \gg  H(P)^{- d +1},
\eqno (10)
$$
valid for any integer polynomial $P(x)$ of degree $d$ having
at least $k \ge 2$
distinct roots $\alpha_1, \ldots , \alpha_k$, established in \cite{Mah64}.

In \cite{BM2} the authors extended the definition of $e(d)$
to clusters of at least $3$ roots.
Let $k$ and $d$ be integers with $2 \le k \le d$.
We denote by $E(d,k)$, respectively
$E_{{\rm irr}}(d,k)$, the infimum of the real numbers $\delta$
for which
$$
\prod_{1\le i<j\le k} \vert \alpha_i -\alpha_j \vert  \ge H(P)^{- \delta}
$$
holds for every integer polynomial $P(x)$, respectively,
irreducible integer polynomial $P(x)$,
of degree $d$ and sufficiently large
height, with distinct roots $\alpha_1, \ldots , \alpha_d$.
We further use the notation
$E^*(d,k)$, respectively $E_{{\rm irr}}^*(d,k)$,
when we restrict our attention to monic
integer polynomials, respectively monic
integer irreducible polynomials.

We deduce from (10) that the quantities $E(d,k)$ and $E^* (d, k)$
are all bounded from above by $d-1$. Regarding lower bounds,
it is proved in \cite{BM1,BM2} that, for any integer $d \ge 4$
and any integer $k\ge 2$ that divides $d$, we have
$$
E(d, k) \ge \frac{k-1}{k} d  \quad \mbox{and} \quad
E^* (d, k) \ge \frac{k-1}{k} d - \frac{k-1}{2}.  \eqno (11)
$$
In fact, the results from \cite{BM1,BM2}
show that these bounds are valid also for $E_{{\rm irr}} (d, k)$ and
$E_{{\rm irr}}^* (d, k)$, respectively.
The constructions presented in the Sections 2 and 3 allow
us to strengthen the estimates (11).

\begin{theorem}
For every integer $k \ge 3$, there exist rational numbers $c_1(k)$ and $c_2(k)$
such that, for every integer $d \ge k$, we have
$$
E(d,k) \geq  \frac{k}{k+1} d - c_1 (k)    \eqno (12)
$$
and
$$
E^* (d,k) \geq  \frac{k}{k+1} d - c_2 (k). \eqno (13)
$$
In particular, the exponent $-d+1$ in (10) cannot be replaced by $- \alpha_k d$,
for a real number $\alpha_k$ less than $k/(k+1)$, even if the
polynomial $P(x)$ in (10) is assumed
to be monic.
\end{theorem}

\noindent{\it Proof. }
Let $\delta \ge 2$ and $h \ge 0$ be integers. The polynomial
$$
{\tilde P}_{\delta, h, n} (x) :=
P_{\delta, n} (x) \times p_{\delta,n} (x) \times \cdots \times p_{\delta+h, n} (x)
$$
has degree $(h+2)\delta + h(h+1)/2$, height $\ll  n^{h+4}$, and
it has a cluster of $h+3$ roots close
to each other. Setting $k = h+3$, a short calculation gives (12), since one can
multiply ${\tilde P}_{\delta, h, n} (x)$ by a suitable power of the monomial $x$.
To get (13), it suffices, for $\delta$ even, to consider the polynomial
$$
{\tilde Q}_{\delta, h, n} (x) :=
Q_{\delta, n} (x) \times q_{\delta,n} (x) \times \cdots \times q_{\delta + 2 h, n} (x),
$$
multiplied by a suitable power of the monomial $x$.
We omit the details.
\qed

\bigskip

{\small \noindent
Yann Bugeaud \\
Universit\'e de Strasbourg  \\
D\'epartement de Math\'ematiques \\ 7, rue Ren\'e Descartes \\
67084 Strasbourg, France \\
{\em E-mail address}: {\tt bugeaud@math.unistra.fr}}

\bigskip

{\small \noindent
Andrej Dujella \\
Department of Mathematics \\ University of
Zagreb
\\ Bijeni\v{c}ka cesta 30 \\
10000 Zagreb, Croatia \\
{\em E-mail address}: {\tt duje@math.hr}}

\end{document}